\documentclass[a4paper]{amsart}
\usepackage{amsmath, tabularx}
\usepackage{amsthm}
\usepackage{amssymb}
\usepackage{amsfonts}
\usepackage{paralist}
\usepackage{aliascnt}
\usepackage{amsrefs}
\usepackage{amscd}
\usepackage{blkarray}

\usepackage{setspace}
\usepackage[colorlinks=true,linkcolor=blue,urlcolor=blue]{hyperref}
\usepackage{tikz}
\usetikzlibrary{matrix}

\newcommand{\RR}{\normalfont\mathbb{F}}
\newcommand{\kk}{\normalfont\mathbb{K}}

\newcommand{\mm}{\normalfont\mathfrak{m}}
\newcommand{\rank}{\normalfont\text{rank}}

\newcommand{\Tor}{\normalfont\text{Tor}}
\newcommand{\Ext}{\normalfont\text{Ext}}

\newcommand{\Ker}{\normalfont\text{Ker}}

\newcommand{\HT}{\normalfont\text{ht}}

\newcommand{\Sym}{\normalfont\text{Sym}}
\newcommand{\Rees}{\normalfont\text{Rees}}
\newcommand{\Hom}{\normalfont\text{Hom}}

\newcommand{\EEQ}{\mathcal{K}}
\newcommand{\REQ}{\mathcal{I}}
\newcommand{\T}{\mathcal{T}}
\newcommand{\M}{\mathcal{M}}
\newcommand{\SE}{\mathcal{S}}
\newcommand{\Fou}{\normalfont\mathcal{F}}
\newcommand{\HL}{\normalfont\text{H}_{\mm}}
\newcommand{\Hdr}{\normalfont\text{H}_{\text{dR}}}
\newcommand{\gr}{\normalfont\text{gr}}
\newcommand{\Sol}{\normalfont\text{Sol}}
\newcommand{\TwistSol}{{\Sol(L_1,L_2;S)}_{\Fou}}
\newcommand{\AAA}{\mathcal{A}}
\newcommand{\BBB}{\mathcal{B}}
\newcommand{\iniTerm}{\normalfont\text{in}}
\newcommand{\bideg}{\normalfont\text{bideg}}
\newcommand{\LCAT}{\M_U^l(\T)}
\newcommand{\RCAT}{\M_U^r(\T)}
\newtheorem{theorem}{Theorem}[section]

\newtheorem{data}{Data}
\newtheorem{headthm}{Theorem}
\newaliascnt{corollary}{theorem}
\newtheorem{corollary}[corollary]{Corollary}
\aliascntresetthe{corollary}

\newaliascnt{lemma}{theorem}
\newtheorem{lemma}[lemma]{Lemma}
\aliascntresetthe{lemma}

\newaliascnt{conjecture}{theorem}

\aliascntresetthe{conjecture}

\newaliascnt{proposition}{theorem}
\newtheorem{proposition}[proposition]{Proposition}
\aliascntresetthe{proposition}

\newaliascnt{definition}{theorem}
\newtheorem{definition}[definition]{Definition}
\aliascntresetthe{definition}

\newaliascnt{notation}{theorem}
\newtheorem{notation}[notation]{Notation}
\aliascntresetthe{notation}

\newaliascnt{example}{theorem}
\newtheorem{example}[example]{Example}
\aliascntresetthe{example}

\newaliascnt{remark}{theorem}
\newtheorem{remark}[remark]{Remark}
\aliascntresetthe{remark}

\newaliascnt{problem}{theorem}

\aliascntresetthe{problem}

\newaliascnt{construction}{theorem}

\aliascntresetthe{construction}

\newaliascnt{defprop}{theorem}
\newtheorem{defprop}[defprop]{Definition-Proposition}
\aliascntresetthe{defprop}

\def\equationautorefname~#1\null{(#1)\null}
\def\sectionautorefname~#1\null{Appendix\null}

\begin{document}

% ---------------Tittle and presentation----------------------------
\title{A $D$-module approach on the equations of the Rees algebra}
\author{Yairon Cid-Ruiz}
\address{Department de Matem\`{a}tiques i Inform\`{a}tica, Facultat de Matem\`{a}tiques i Inform\`{a}tica, Universitat de Barcelona, Gran Via de les Corts Catalanes, 585; 08007 Barcelona, Spain.}
\email{ycid@ub.edu}
\urladdr{http://www.ub.edu/arcades/ycid.html}
\date{\today}



\thanks{The author was funded by the European Union's Horizon 2020 research and innovation programme under the Marie Sk\l{}odowska-Curie grant agreement No. 675789.}

\begin{abstract}
	Let $I \subset R = \RR[x_1,x_2]$ be a height two ideal minimally generated by three homogeneous polynomials of the same degree $d$, where $\RR$ is a field of characteristic zero.
	We use the theory of $D$-modules to deduce information about the defining equations of the Rees algebra of $I$.
	Let $\EEQ$ be the kernel of the canonical map $\alpha: \Sym(I) \rightarrow \Rees(I)$ from the symmetric algebra of $I$ onto the Rees algebra of $I$.
	We prove that $\EEQ$ can be described as the solution set of a system of differential equations, that the whole bigraded structure of $\EEQ$ is characterized by the integral roots of certain $b$-functions, and that certain de Rham cohomology groups can give partial information about $\EEQ$. 
\end{abstract}
\maketitle

\section{Introduction}
Let $\RR$ be a field of characteristic zero, $R=\RR[x_1,x_2]$ a polynomial ring in two variables, and $I=(f_1,f_2,f_3) \subset R$ be a height two ideal minimally generated by three homogeneous polynomials of the same degree $d$.
The Rees algebra of $I$ is defined as $\Rees(I)=R[It]=\bigoplus_{i=0}^{\infty}I^it^i$. We can see $\Rees(I)$ as a quotient of the polynomial ring $S=R[T_1,T_2,T_3]$ via the map
\begin{equation}
	\label{Rees_presentation}
	S=R[T_1,T_2,T_3] \xrightarrow{\psi} \Rees(I), \quad \psi(T_i)=f_it.
\end{equation}
Of particular interest are the defining equations of the Rees algebra $\Rees(I)$, that is, the kernel $\REQ=\Ker(\psi)$ of this map $\psi$. 
A large number of works have been done to determine the equations of the Rees algebra, and the problem has been studied by algebraic geometers and commutative algebraists under various conditions (see e.g. \cite{VASCONCELOS_BLOWUP} and the references therein).
In recent years, a lot of attention has been given to find the minimal generators of the equations of the Rees algebra for an ideal in a polynomial ring (see e.g. \cite{LAURENT, CARLOS_CONICS, CARLOS_MU_TWO, MONOMIAL_CARLOS, COX_MOVING, COX_HOFFMAN_WANG,KPU_BIGRAD, KPU_SCROLL,KPU_Gor_3, LIN_POLINI, MADSEN,VASC_SIMIS_HONG}), partly inspired by new connections with geometric modelling.
Despite this extensive effort, even in the ``simple'' case studied in this paper, the problem of finding the minimal generators of $\REQ$ remains open.

By the Hilbert-Burch theorem we know that the presentation of $I$ is of the form
\begin{equation}
\label{HilbertBurchRes}
0 \rightarrow R(-d-\mu)\oplus R(-2d+\mu) \xrightarrow{\varphi} R(-d)^3 \xrightarrow{\left[f_1,f_2,f_3\right]} I \rightarrow 0,
\end{equation}
and $I$ is generated by the $2 \times 2$-minors of $\varphi$; we may assume that $0 < \mu \le d-\mu$.
The symmetric algebra of $I$ can easily be described as a quotient of $S$ by using the presentation of $I$ . 
We define the equations of the symmetric algebra as
\begin{equation}
\label{linear_eqs}
 \left[g_1, g_2\right] = \left[T_1, T_2, T_3\right]\cdot \varphi, 
\end{equation}
then we get $\Sym(I) \cong S/(g_1,g_2)$. There is an important relation between $\Sym(I)$ and $\Rees(I)$ in the form of the following canonical exact sequence 
$$
0 \rightarrow \EEQ \rightarrow \Sym(I) \xrightarrow{\alpha} \Rees(I) \rightarrow 0.
$$
Here we have $\EEQ=\REQ/(g_1,g_2)$, which allows us to take $\EEQ$ as the object of study.

We introduce the polynomial ring $U=\RR[T_1,T_2,T_3]$ and we give a bigraded structure to $S=R \otimes_{\RR} U$, where $\bideg(T_i)=(1,0)$ and $\bideg(x_i)=(0, 1)$; then we denote by $S_{p,q}$ the $\RR$-vector space spanned by the monomials  $x_1^{\alpha_1}x_2^{\alpha_2}T_1^{\gamma_1}T_2^{\gamma_2}T_3^{\gamma_3}$ with $\gamma_1+\gamma_2+\gamma_3=p$ and $\alpha_1+\alpha_2=q$.  
The map $\psi$ from \autoref{Rees_presentation} becomes bihomogeneous when we declare $\bideg(t)=(1,-d)$, and also from the fact that $\bideg(g_1)=(1, \mu)$ and $\bideg(g_2)=(1, d-\mu)$, then we get that $\Rees(I)$, $\Sym(I)$, $\REQ$ and $\EEQ$ have natural structures as bigraded $S$-modules.
For an arbitrary bigraded $S$-module $N$ we use the notations
$$
	N_{p,*}=\bigoplus_{q \in \mathbb{Z}} N_{p,q} \qquad\text{ and }\qquad N_{*,q}=\bigoplus_{p\in \mathbb{Z}} N_{p,q},
$$
where $N_{p,*}$ is a graded $R$-module  and $N_{*,q}$ is a graded $U$-module. 

The main feature of this paper is the use of the theory of $D$-modules in the problem of finding the equations of $\Rees(I)$.
For this we need to introduce the Weyl algebra $A_2(\RR)$ (\autoref{Weyl_algebra}) and  a polynomial ring $\T=A_2(\RR)[T_1,T_2,T_3]$ (\autoref{polynomial_Weyl}) over the Weyl algebra; also we define two differential operators $L_1=\Fou(g_1)$ and  $L_2=\Fou(g_2)$ by applying the Fourier transform (\autoref{Fourier_transform}) to $g_1$ and $g_2$ from \autoref{linear_eqs}.

Our first main result claims that $\EEQ$ can be described by solving a system of differential equations.
\begin{headthm}
	\label{thmA}
	(\autoref{isom_Sol})
	Let $I \subset R=\RR[x_1,x_2]$ be a height two ideal minimally generated by three homogeneous polynomials of the same degree $d$, and let $L_1=\Fou(g_1)$ and $L_2=\Fou(g_2)$ be the Fourier transform of $g_1$ and $g_2$ from \autoref{linear_eqs}. Then we have the following isomorphism of bigraded $S$-modules 
	$$
	\EEQ \cong {\Sol\Big(L_1,L_2;S\Big)}_{\Fou}(-2,-d+2), 
	$$ 
	where ${\Sol(L_1,L_2;S)}=\{ h \in S \mid L_1 \bullet h = 0 \text{ and } L_2 \bullet h = 0\}$ and the subscript-$\Fou$ denotes the twisting by the Fourier transform (\autoref{BIGRAD_STRUCT_SOL}).
\end{headthm}

Since $g_1$ and $g_2$ generate all the linear part of $\REQ$ (the syzygies of $I$) and $\EEQ=\REQ/(g_1,g_2)$, then we have $\EEQ_{p,*}=0$ for all $p < 2$.
As an application of \autoref{thmA} we give a complete characterization of  the graded structure of each $R$-module $\EEQ_{p,*}$ ($p\ge2$) in terms of the integral roots of certain $b$-functions (\autoref{bFunction_module}). 

\begin{headthm}
	\label{thmB}
	(\autoref{expression_bFunction})
	Let $I \subset R=\RR[x_1,x_2]$ be as in \autoref{thmA}.  Then for each integer $p \ge 2$ there exists a nonzero $b$-function $b_p(s)$, and we have a relation between the graded structure of $\EEQ_{p,*}$ and the integral roots of $b_p(s)$ given in the following equivalence
	$$
	\EEQ_{p,q} \neq 0 \qquad \Longleftrightarrow \qquad b_p(-d+2+q)=0.
	$$
	Even more, we have that these are the only possible roots of $b_p(s)$, that is
	$$
	b_p(s)=\prod_{\{q \in \mathbb{Z}\;\mid\; \EEQ_{p,q}\neq 0\}} (s+d-2-q).
	$$
\end{headthm}
This \autoref{thmB}  is interesting for us in the sense that gives a tool for deducing information about $\EEQ$, but on the other hand, from a $D$-module point of view it is worthy to note that describes the $b$-function of a family of holonomic $D$-modules like those in \autoref{notations_groebner}. 

In our last main result we change the role of $L_1$ and $L_2$, more specifically, instead of having them as operators we place them in a quotient. 
We make this change by means of a duality proven in \autoref{duality}, and it allows us to establish an isomorphism of graded $U$-modules between $\EEQ$ and a certain de Rham cohomology group. 
In particular, this isomorphism could give an alternative way to compute or estimate the dimension $\dim_{\RR}(\EEQ_{p,*})$ of each $\EEQ_{p,*}$ regarded as finite dimensional $\RR$-vector space (see \autoref{bigrad_struct_EEQ}).

\begin{headthm}
	\label{thmC}
	(\autoref{deRham_cohomology})
	Let $I\subset R=\RR[x_1,x_2]$, $L_1$ and $L_2$ be as in \autoref{thmA}, and let $Q$ be the left $\T$-module $Q=\T/\T(L_1,L_2)$. 
	Then we have the following isomorphism of graded $U$-modules
	$$\EEQ \cong \Hdr^0(Q) = \{w \in Q \mid \partial_1 \bullet w= 0 \text{ and } \partial_2 \bullet w= 0 \}.
	$$
	In particular, for any integer $p$ we have an isomorphism of $\RR$-vector spaces
	$$
	\EEQ_{p,*} \cong \Hdr^0(Q_p) = \{w \in Q_p \mid \partial_1 \bullet w= 0 \text{ and } \partial_2 \bullet w= 0 \}.
	$$
\end{headthm}

Our next theorem works over an arbitrary field $\kk$. 
We remark that this result it is already known (see e.g. \cite[Lemma 2.4]{COX_HOFFMAN_WANG} or \cite[Theorem 2.4]{KPU_BIGRAD}), but we present a different proof. The rest of our work will depend on it.
\begin{headthm}
	\label{thmD}
	(\autoref{eqs_whole})
	Let $I \subset R=\kk[x_1,x_2]$ be a height two ideal minimally generated by three homogeneous polynomials of the same degree $d$, and let $g_1$ and $g_2$ be as in \autoref{linear_eqs}.
	Then we have the following isomorphism of bigraded $S$-modules
	\begin{equation*}
	\EEQ \cong \Big\{ w \in \HL^2\big(S\big)(-2, -d) \;\big| \; g_1\cdot w =0 \text{ and } g_2\cdot w=0 \Big\}.
	\end{equation*}	
\end{headthm}

The basic outline of this paper is as follows. In \hyperref[section2]{Section 2} we give a proof of \autoref{thmD}.
In \hyperref[section3]{Section 3}
we make a translation of our problem into the theory of $D$-modules and we prove \autoref{thmA}.
In \hyperref[section4]{Section 4} we prove \autoref{thmB}.
In \hyperref[section5]{Section 5} we prove \autoref{thmC}.
In \hyperref[section6]{Section 6} we present a script in \textit{Macaulay2} \cite{M2} that can compute each $b$-function $b_p(s)$ from \autoref{thmB}, and using it we effectively recover the bigraded structure of $\EEQ$ for a couple of examples.

\section{An ``explicit'' description of the equations}\label{section2}
In this section we shall use the following data.
\begin{data}
	\label{first_data}
	Let $\kk$ be an arbitrary field, and $R=\kk[x_1,x_2]$ be the polynomial ring in two variables. 
	Let $I\subset R$ be a height two ideal minimally generated by three homogeneous polynomials $\{f_1,f_2,f_3\}$ of the same degree $d$.
	From Hilbert-Burch theorem we have a presentation
	$$
	0 \rightarrow R(-d-\mu)\oplus R(-2d+\mu) \xrightarrow{\varphi} R(-d)^3 \rightarrow I \rightarrow 0,
	$$
	where the elements of the first column of $\varphi$ are homogeneous of degree $\mu$, and the elements of the second column are homogeneous of degree $d-\mu$.
	Let $U$ and $S$ be the polynomial rings 
	$U=\kk[T_1,T_2,T_3]$ and $S=R[T_1,T_2,T_3]=\kk[x_1,x_2,T_1,T_2,T_3]$. 
	We regard $S$ as a bigraded $\kk$-algebra, where $\bideg(T_i)=(1, 0)$ and $\bideg(x_i)=(0, 1)$. 
	The equations of the symmetric algebra are given by $\left[g_1,g_2\right] = \left[T_1,T_2,T_3\right]\cdot \varphi.$
	We are interested in the kernel $\EEQ$ of the surjective map $	\alpha: \Sym(I) \rightarrow \Rees(I).$	
\end{data}

By the pioneering work of \cite{MICALI_REES}, we can compute $\EEQ$ as the torsion in $\Sym(I)$ with respect to the maximal ideal $\mm=(x_1,x_2) \subset R$, that is 
$$
\EEQ = (0 :_{\Sym(I)} \mm^{\infty}) = \HL^0(\Sym(I)).
$$

Given a bigraded $S$-module $M$, by definition each  local cohomology module $\HL^j(M)$ is only an $R$-module.
In the following lemma we endow $\HL^j(M)$ with a structure of bigraded $S$-module. We use \cite[Chapter 13]{Brodmann_Sharp_local_cohom}
for the foundations of local cohomology modules in the graded case.

\begin{lemma}
	Let $M$ be a bigraded $S$-module. Then, the following statements hold:
	\begin{enumerate}[(i)]
		\item Use the decomposition $M = \bigoplus_{p\in \mathbb{Z}} M_p$, where $M_p$ is the graded $R$-module given by $M_p = \bigoplus_{q \in \mathbb{Z}} M_{p,q}$. Then 
		\begin{equation}
				\label{local_cohom_bigrad}
			\HL^j(M)=\bigoplus_{p \in \mathbb{Z}} \HL^j(M_p), 			
		\end{equation}	
		is a bigraded $S$-module with ${\HL^j(M)}_{p,q}={\HL^j(M_p)}_q$ (where ${\HL^j(M_p)}_q$ represents the $q$-th graded part of the graded $R$-module $\HL^j(M_p)$).
		The actions of the $x_i$'s are natural because $\HL^j(M)$ is an $R$-module. 
		The action of the $T_i$'s over $M$ can be seen as homogeneous homomorphisms $T_i: M_p \rightarrow M_{p+1}$ of graded $R$-modules, then the induced homogeneous homomorphisms $\HL^j(T_i): \HL^j(M_p) \rightarrow \HL^j(M_{p+1})$  of graded $R$-modules give us the action of the $T_i$'s over $\HL^j(M)$.
		\item For any $a, b \in \mathbb{Z}$, we have the isomorphism of bigraded $S$-modules $\HL^j(M(a,b))\cong\HL^j(M)(a,b)$.
 	\end{enumerate}
 	\begin{proof}
		$(i)$ The decomposition \autoref{local_cohom_bigrad} comes from the fact that local cohomology commutes with direct sums, and that each $\HL^j(M_p)$ has a natural structure of graded $R$-module (see e.g. \cite[Chapter 13]{Brodmann_Sharp_local_cohom}).
		
		$(ii)$ The shifting on the $T_i$'s follows from the construction \autoref{local_cohom_bigrad} and so we are left to check that $\HL^j(M_p(b)) \cong \HL^j(M_p)(b)$ for each $p \in \mathbb{Z}$. 
		For this, we use \cite[Theorem 13.4.5]{Brodmann_Sharp_local_cohom} and any of the remarks in page 273 of \cite{Brodmann_Sharp_local_cohom}, for instance using the construction as a direct limit of Ext's we have
		$$
		\HL^j(M_p(b)) \cong \varinjlim_{n} {}^{*}\Ext_R^j(R/\mm^n, M_p(b))
		\cong  \varinjlim_{n} {}^{*}\Ext_R^j(R/\mm^n, M_p)(b) \cong \HL^j(M_p)(b).
		$$
		(see e.g. \cite[Section 1.5]{BRUNS_HERZ} for graded dual ${}^{*}\Hom_R$ and its derived functors ${}^{*}\Ext_R^j$ in the category of graded modules).
 	\end{proof} 
\end{lemma}

The ``philosophy'' that we follow in this section is similar to the one used in \cite{MADSEN}. Explicitly, we shall try to find information by deleting the columns of $\varphi$ and hopefully work with ``simpler'' modules. 
Let $\varphi_1$ be the matrix given by the first column of $\varphi$, then we we are interested in the module $E = \text{Coker}(\varphi_1)$ with presentation 
$$
0 \rightarrow R(-d-\mu) \xrightarrow{\varphi_1} R(-d)^3 \rightarrow E \rightarrow 0.  
$$
\begin{lemma}
	\label{describe_E}
	For the module $E$ we have
	\begin{enumerate}[(i)]
		\item	$\Sym(E) \cong S/(g_1)$;
		\item $\Sym(E)$ is an integral domain.
	\end{enumerate}
	\begin{proof}
		$(i)$ Follows from the  presentation of $E$.
		
		$(ii)$ Since $I_1(\varphi_1) \supset I_2(\varphi)$, we have that $\HT(I_1(\varphi_1))=2$. Then, by
		\cite[Theorem 3.4]{SIMIS_VASC_SYM_INT} we get that $\Sym(E)$ is an integral domain.
	\end{proof}
\end{lemma}

Now we can find explicit relations between the local cohomology modules of $\Sym(I)$ and $\Sym(E)$ from the important fact that $\Sym(E)$ is an integral domain.

\begin{lemma}
	We have the following exact sequences of bigraded $S$-modules
	\begin{align}
			&\label{exact_seq_Sym_p_I}
			0 \rightarrow \HL^0(\Sym(I)) \xrightarrow{\partial} \HL^1\big(\Sym(E)\big)(-1,-d+\mu) \xrightarrow{g_2} \HL^1(\Sym(E));\\
			&\label{exact_seq_Sym_p_E}
			0 \rightarrow \HL^1(\Sym(E)) \xrightarrow{\partial} \HL^2\big(S\big)(-1, -\mu) \xrightarrow{g_1} \HL^2(S).
		\end{align}
	\begin{proof}
		Since $\Sym(E)\cong S/(g_1)$ is an integral domain we have a short  exact sequence 
		$$
		0 \rightarrow \Sym\big(E\big)(-1, -d+\mu) \xrightarrow{g_2} \Sym(E) \rightarrow \Sym(I) \rightarrow 0.
		$$
		Using the corresponding long exact sequence in local cohomology and the fact that $\HL^0(\Sym(E))=0$, we get the required exact sequence 
		$$
		0 \rightarrow \HL^0(\Sym(I)) \xrightarrow{\partial} \HL^1\big(\Sym(E)\big)(-1,-d+\mu) \xrightarrow{g_2} \HL^1(\Sym(E)),
		$$
		where $\partial$ is the induced connecting homomorphism.
		
		Similarly, from the short exact sequence 
		\begin{equation}
			\label{res_SymE}
			0 \rightarrow S(-1,-\mu) \xrightarrow{g_1} S \rightarrow \Sym(E) \rightarrow 0,
		\end{equation} 
		and the fact that 
		$$
		\HL^j(R)\cong \begin{cases}
		x_1^{-1}x_2^{-1}\kk[x_1^{-1},x_2^{-1}]  \qquad \text{if } j = 2\\
		0 \qquad \qquad \qquad\qquad\quad\text{otherwise,}
		\end{cases}
		$$
		we can follow the same long exact sequence argument and obtain \autoref{exact_seq_Sym_p_E}.
	\end{proof}
\end{lemma}

The next theorem contains the main result of this section, where we find an ``explicit'' way of computing the equations of the Rees algebra of $I$.

\begin{theorem}
	\label{eqs_whole}
	Adopt \autoref{first_data}. 
	Then we have the following isomorphism of bigraded $S$-modules
	\begin{equation*}
		\EEQ \cong \Big\{ w \in \HL^2\big(S\big)(-2, -d) \;\big| \; g_1\cdot w =0 \text{ and } g_2\cdot w=0 \Big\}.
	\end{equation*}
	\begin{proof}
		The  commutative diagram
		\begin{center}		
			\begin{tikzpicture}
			\matrix (m) [matrix of math nodes,row sep=2em,column sep=4em,minimum width=2em, text height=1.5ex, text depth=0.25ex]
			{
				S(-2, -d) & S(-1, -d+\mu)  \\
				S(-1, -\mu) & S \\
			};
			\path[-stealth]
			(m-1-1) edge node [above] {$g_1$} (m-1-2)
			(m-2-1) edge node [above] {$g_1$} (m-2-2)
			(m-1-1) edge node [right]	{$g_2$} (m-2-1)
			(m-1-2) edge node [right]	{$g_2$} (m-2-2)
			;				
			\end{tikzpicture}	
		\end{center}
		can be extended to the following commutative diagram with exact rows (each row is as in \autoref{res_SymE})
	\begin{center}	
		\begin{tikzpicture}
		\matrix (m) [matrix of math nodes,row sep=2em,column sep=2em,minimum width=2em, text height=1.5ex, text depth=0.25ex]
		{
			0 & S(-2,-d) & S(-1,-d+\mu) & \Sym\big(E\big)(-1,-d+\mu) & 0\\
			0 & S(-1,-\mu) & S & \Sym(E) & 0. \\
		};
		\path[-stealth]
		(m-1-1) edge (m-1-2)
		(m-2-1) edge (m-2-2)
		(m-1-2) edge node [above] {$g_1$} (m-1-3)
		(m-2-2) edge node [above] {$g_1$} (m-2-3)
		(m-1-2) edge node [right]	{$g_2$} (m-2-2)
		(m-1-3) edge node [right]	{$g_2$} (m-2-3)
		(m-1-4) edge node [right]	{$g_2$} (m-2-4)
		(m-1-3) edge (m-1-4)
		(m-2-3) edge (m-2-4)
		(m-1-4) edge (m-1-5)
		(m-2-4) edge (m-2-5)
		;				
		\end{tikzpicture}	
	\end{center}
	From the ``naturality of the connecting homomorphism $\partial$'' \cite[Chapter 6]{ROTMAN} and \autoref{exact_seq_Sym_p_E}, we get the following commutative diagram with exact rows
		\begin{center}		
		\begin{tikzpicture}
		\matrix (m) [matrix of math nodes,row sep=2em,column sep=2em,minimum width=2em, text height=1.5ex, text depth=0.25ex]
		{
			0 & \HL^1\big(\Sym(E)\big)(-1,-d+\mu) & \HL^2\big(S\big)(-2,-d) & \HL^2\big(S\big)(-1,-d+\mu)\\
			0 & \HL^1(\Sym(E))  & \HL^2\big(S\big)(-1,-\mu) & \HL^2(S).\\
		};
		\path[-stealth]
		(m-1-1) edge (m-1-2)
		(m-2-1) edge (m-2-2)
		(m-1-3) edge node [above] {$g_1$} (m-1-4)
		(m-2-3) edge node [above] {$g_1$} (m-2-4)
		(m-1-2) edge node [right]	{$g_2$} (m-2-2)
		(m-1-3) edge node [right]	{$g_2$} (m-2-3)
		(m-1-4) edge node [right]	{$g_2$} (m-2-4)
		(m-1-2) edge node [above] {$\partial$} (m-1-3)
		(m-2-2) edge node [above] {$\partial$} (m-2-3)
		;				
		\end{tikzpicture}	
	\end{center}
	From this diagram and \autoref{exact_seq_Sym_p_I}, we get the exact sequence 
	$$
	0 \rightarrow \EEQ \rightarrow \Ker\Big(\HL^2\big(S\big)(-2,-d) \xrightarrow{g_2} \HL^2\big(S\big)(-1,-\mu)\Big) \xrightarrow{g_1} \Ker\Big(\HL^2\big(S\big)(-1,-d+\mu) \xrightarrow{g_2} \HL^2(S)\Big), 
	$$
	from which we finally identify
	$$
	\EEQ \cong \Big\{ w \in \HL^2\big(S\big)(-2, -d) \;\big| \; g_1\cdot w =0 \text{ and } g_2 \cdot w=0 \Big\}.
	$$
	\end{proof}
\end{theorem}

\begin{corollary}
	\label{bigrad_struct_EEQ}
	Adopt \autoref{first_data}.
	The following statements hold:
	\begin{enumerate}[(i)]
		\item For $p \ge 2$ the graded part $\EEQ_{p,*}$ is a finite dimensional $\kk$-vector space with $\EEQ_{p,d-2} \neq 0$ and $\EEQ_{p,q}=0$ for $q > d-2$.
		\item $\EEQ_{*,d-2}\cong U(-2)$  is an isomorphism of graded $U$-modules.
	\end{enumerate}
	\begin{proof}						
		$(i)$ For any $q > d - 2$ we have  $q-d>-2$, and so ${\HL^2(R)}_{q-d}=0$ which implies $\EEQ_{p,q}=0$.
		If $q = d-2$ then ${\HL^2(R)}_{q-d}=\kk\cdot\frac{1}{x_1x_2}$, and so it follows that $\EEQ_{p,q}\neq 0$ since $x_1 \cdot \frac{1}{x_1x_2}=0$ and $x_2\cdot \frac{1}{x_1x_2}=0$.
		
		$(ii)$ It follows from the fact that $S\cdot\frac{1}{x_1x_2} \cong \kk[T_1,T_2,T_3]=U$.	
	\end{proof}
\end{corollary}

In this previous corollary we have seen that the maximal $x$-degree of every graded part $\EEQ_{p,*}$ is the same and equal to $d-2$, 
but for the minimal $x$-degree of $\EEQ_{p,*}$ there is no such nice characterization.
In \hyperref[section4]{Section 4} under the assumption of working over a field of characteristic zero, we shall relate the minimal $x$-degree with the integral roots of certain $b$-functions.

\section{Translation into $D$-modules}\label{section3}

The core of this section is to translate our problem into $D$-modules. 
A good introduction to the theory of $D$-modules can be found in  \cite{BJORK} or \cite{COUTINHO}.
The section is divided into two subsections, a first one containing some notations and definitions regarding $D$-modules that we shall use for  the rest of this paper, and a second one containing our translation. 

\subsection{Notations}\label{notations_sec}
For the rest of this paper we shall work over a field $\RR$ of characteristic zero, and from now on we shall use the following data.

\begin{data}
	\label{data}
	Adopt \autoref{first_data} and change the arbitrary field $\kk$ for a field $\RR$ of characteristic zero.
\end{data} 

We introduce the ring of $\RR$-linear differential operators over $R=\RR[x_1,x_2]$, which in our characteristic zero case coincides with the Weyl algebra.

\begin{definition}
	\label{Weyl_algebra}
	The Weyl algebra $D=A_2(\RR)$ is defined as a quotient of the free algebra $\RR \textless x_1,x_2,\partial_1,\partial_2\textgreater$ by the two sided ideal generated by the relations
	$$
	x_ix_j=x_jx_i, \qquad \partial_i \partial_j = \partial_j \partial_i, \qquad \partial_ix_j = x_j\partial_i + \delta_{ij},
	$$
	where $\delta_{ij}$ is Kronecker's symbol.
\end{definition}

The $D$-module structure of $R$ is given by: for any $f \in R$, the operator $x_i$ is the usual multiplication  $x_i\bullet f = x_if$ and the operator $\partial_i$ is the differentiation $\partial_i\bullet f=\frac{\partial f}{\partial x_i}$. We shall always stress the action of the Weyl algebra by using the symbol ``$\;\bullet$''.
Thus, for instance, if we regard $x_1 \in R$ then we have $\partial_1 \bullet x_1=1$, but instead for $x_1 \in D$ we have $\partial_1x_1=x_1\partial_1 + 1$.

Of particular interest are the holonomic $D$-modules. A finitely generated left $D$-module $M\neq 0$ is said to be holonomic if it has Bernstein dimension $d(M)=2$, or equivalently, if $\Ext_D^i(M, D)$ vanishes for all $i \neq 2$.
A left $D$-ideal $J$ is said to be holonomic when $D/J$ is holonomic.

All the modules in the \v{C}ech complex are localizations of $R$, thus by defining the $D$-module structure of any localization $R_f$ of $R$, the local cohomology modules obtain a natural structure as $D$-modules (see e.g. \cite[Lecture 23]{24H_LOC_COHOM}). 
For any localization $R_f$ the $D$-module structure is defined by
$$
x_i \bullet \frac{g}{f^k}=x_i\frac{g}{f^k} \qquad \text{and} \qquad \partial_i \bullet \frac{g}{f^k} = \frac{1}{f^k}\frac{\partial g}{\partial x_i} - \frac{kg}{f^{k+1}}\frac{\partial f}{\partial x_i}.
$$

Due to the non-commutativity of $D$, we need to take some care with the maps of left or right $D$-modules. 
Let $A \in D^{r\times s}$ be an $r\times s$ matrix with entries in $D$.
Multiplying with $A$ gives us a map of left $D$-modules,
$$
D^r \xrightarrow{\cdot A} D^s \quad : \quad \left[\ell_1,\ldots,\ell_r \right] \mapsto  \left[\ell_1,\ldots,\ell_r \right] \cdot A,
$$
where we regard $D^r$ and $D^s$ as row vectors.

The matrix $A \in D^{r \times s}$ also defines a map of right $D$-modules in the opposite direction, 
$$
{(D^s)}^T \xrightarrow{A\cdot} {(D^r)}^T \quad : \quad {\left[\ell_1^{'},\ldots,\ell_s^{'} \right]}^T \mapsto  A \cdot {\left[\ell_1^{'},\ldots,\ell_s^{'} \right]}^T,
$$
where the superscript-$T$ means that ${(D^s)}^T$ and ${(D^r)}^T$ are considered as column vectors. 
The right $D$-module ${(D^s)}^T$ may be regarded as the dual module  $\Hom_D(D^s,D)$.
Applying $\Hom_D(-,D)$ to the map $D^r \xrightarrow{\cdot A} D^s$ of left $D$-modules induces the map ${(D^s)}^T \xrightarrow{A\cdot} {(D^r)}^T$ of right $D$-modules.

We have an equivalence between the category of left $D$-modules and the category of right $D$-modules, given by the algebra involution 
$$
D \xrightarrow{\tau} D \quad : \quad x^{\alpha}\partial^\beta \mapsto {(-\partial)}^\beta x^\alpha.
$$
The map $\tau$ is called the standard transposition. 
For instance, given a left $D$-module $D^r/M_0$ its corresponding standard transposition is the right $D$-module
$$
\tau\Big(\frac{D^r}{M_0}\Big) = \frac{D^r}{\tau(M_0)}, \qquad \tau(M_0)=\{\tau(L) \mid L \in M_0 \}.
$$
See \cite{OAKU_GROBNER} for more details on the standard transposition $\tau$.

Finally, to describe all the graded parts $\EEQ_{p,*}$ together, we need to define a larger algebra to work in. 

\begin{definition}
	\label{polynomial_Weyl}
	We define $\T$ as a polynomial ring in the three variables $T_1,T_2,T_3$ over the Weyl algebra, that is $\T = A_2(\RR)[T_1,T_2,T_3]=\RR[x_1,x_2]\textless\partial_1,\partial_2\textgreater
	[T_1,T_2,T_3]$.
\end{definition}

We extend the standard transposition $\tau$ over $\T$ by making $\tau(T_i)=T_i$. 
This algebra $\T$ is naturally a graded $U$-module with grading on the $T_i$'s, and by $\T_p$ we denote the free $D$-module spanned by the monomials $\mathbf{T}^\gamma$ with $\lvert\gamma \rvert=p$, that is, $\T_p=D^{\binom{p+2}{2}}$.
Also, for technical purposes we shall need to introduce the subcategory $\LCAT$ of left $\T$-modules with an underlying structure of graded $U$-module.
The subcategory $\RCAT$ of $U$-graded right $\T$-modules can be defined in a completely similar way.
We essentially follow the exposition of \cite[Section 1.5]{BRUNS_HERZ}.

\begin{definition}
	We say that a left $\T$-module $M$ has an underlying structure of graded $U$-module (or simply that it is $U$-graded) when it has a decomposition $M = \bigoplus_{i \in \mathbb{Z}} M_i$, where each $M_i$ is a left $D$-module and $\T_p \bullet M_i \subset M_{i+p}$.
\end{definition}

\begin{definition}
	The category $\LCAT$, has as objects the left $\T$-modules with an underlying structure of graded $U$-module.
	A morphism $\varphi: M \rightarrow N$ in $\LCAT$ is a homomorphism of left $\T$-modules satisfying $\varphi(M_i) \subset N_i$ for all $i \in \mathbb{Z}$.
\end{definition}

If $M$ belongs to $\LCAT$, then $M(i) \in \LCAT$ denotes the $U$-graded left $\T$-module with grading given by $M(i)_n=M_{i+n}$.
All the following assertions follow from the fact that the $T_i$'s are central in $\T$.

Since each module $M \in \LCAT$ is a homomorphic image (in $\LCAT$) of a free module (in $\LCAT$) of the form $\bigoplus \T(i)$ (simply by choosing homogeneous generators of $M$), then the category $\LCAT$ has enough projectives.
Thus, every module $M \in \LCAT$ has a free resolution in $\LCAT$, and this fact allows us to define derived functors in $\LCAT$ (see e.g. \cite[Chapter 6]{ROTMAN}).

Let $M \in \RCAT$ be a $U$-graded right $\T$-module and $N \in \LCAT$ be a $U$-graded left $\T$-module. Then from the non-commutativity of $D$ follows that the tensor product $M \otimes_\T N$ has only a structure of graded $U$-module; its homogeneous component $(M\otimes_\T N)_n$ is generated (as an $\RR$-vector space) by the elements $u \otimes_\T v$ with $u \in M_i$, $v \in N_j$ and $i + j=n$.
Using that each module in $\LCAT$ or in $\RCAT$ has a free resolution (in $\LCAT$ or in $\RCAT$), then (the $\RR$-vector space) $\Tor_i^{\T}(M, N)$ has a natural structure of graded $U$-module for any $i \ge 0$.  
We shall use the notation ${}^*\Tor_i^{\T}(M, N)$ to stress its graded structure as a $U$-module.

Let $M, N \in \LCAT$ be $U$-graded left $\T$-modules.  
A homomorphism of left $\T$-modules $\varphi:M \rightarrow N$ is called homogeneous of degree $i$ if $\varphi(M_n) \subset N_{n+i}$ for all $n \in \mathbb{Z}$.
We denote by $\Hom_i(M, N)$ the $\RR$-vector space of homogeneous homomorphisms of degree $i$. 
The $\RR$-vector subspaces $\Hom_i(M, N)$ of $\Hom_\T(M, N)$ form a direct sum, and we have that 
$$
{}^*\Hom_{\T}(M, N)=\bigoplus_{i \in \mathbb{Z}} \Hom_i(M, N)
$$
is naturally a graded $U$-module.
Also, when $M$ is finitely generated we have that ${}^*\Hom_{\T}(M, N)=\Hom_{\T}(M, N)$.

For any $N \in \LCAT$ we define ${}^*\Ext_\T^i(M, N)$ as the $i$-th right derived functor of ${}^*\Hom_{\T}(-, N)$ in $\LCAT$.
Hence, given a projective resolution $P_{\bullet}$ of $M$ in $\LCAT$, we have
$$
{}^*\Ext_\T^i(M, N) \cong H^i({}^*\Hom_\T(P_{\bullet}, N)),
$$
for all $i \ge 0$.
A particular and important case is when $N = \T$, since $\T$ can be seen as a bimodule then we have that ${}^*\Ext_{\T}^i(M, \T)$ is a module in the category $\RCAT$ of right $\T$-modules with a structure of graded $U$-module.

The Weyl algebra $D=A_2(\RR)$ is a left Noetherian ring (see e.g. \cite[Proposition 2.8, page 6]{BJORK}), then from the Hilbert basis theorem (see e.g. \cite[Theorem 3.21]{ROTMAN}) we have that $\T$ is also a left Noetherian ring. 
Thus, for $M \in \LCAT$ finitely generated we can find a resolution in $\LCAT$ made-up of finitely generated free modules,
and so we have that ${}^*\Ext_{\T}^i(M, N)=\Ext_{\T}^i(M, N)$.
We shall use the notation ${}^*\Ext_{\T}^i(M, N)$ to emphasize its graded structure as a $U$-module.

\subsection{The translation}

We can see that $S=\bigoplus_{\gamma} R\mathbf{T}^\gamma$ and $\HL^2(S)=\bigoplus_{\gamma}\HL^2(R)\mathbf{T}^\gamma$
both belong to the category $\LCAT$ of $U$-graded left $\T$-modules.

\begin{proposition}
	\begin{enumerate}[(i)]
		\item The left $\T$-module $\HL^2(S)$ is cyclic with generator $\frac{1}{x_1x_2}$ and presentation 
		$$
		0 \rightarrow \T(x_1,x_2) \rightarrow \T \xrightarrow{\bullet \frac{1}{x_1x_2}} \HL^2(S) \rightarrow 0.
		$$
		\item The left $\T$-module $S$ is cyclic with generator $1$ and presentation 
		$$
		0 \rightarrow \T(\partial_1,\partial_2) \rightarrow \T \xrightarrow{\bullet 1} S \rightarrow 0.
		$$	
	\end{enumerate}
	\begin{proof}
	$(i)$ To prove that $\frac{1}{x_1x_2}$ is a generator of $\HL^2(S)$ it is enough to show that any monomial $\frac{1}{x_1^{\alpha_1}x_2^{\alpha_2}}T_1^{\gamma_1}T_2^{\gamma_2}T_3^{\gamma_3}$ belongs to $\T \bullet \frac{1}{x_1x_2}$, but this is obtained from the fact that $\text{char}(\RR)=0$ and the following identity
	$$
	\partial_1^{\alpha_1-1}\partial_2^{\alpha_2-1}T_1^{\gamma_1}T_2^{\gamma_2}T_3^{\gamma_3} \bullet \frac{1}{x_1x_2}=(-1)^{\alpha_1+\alpha_2}\frac{(\alpha_1-1)!(\alpha_2-1)!}{x_1^{\alpha_1}x_2^{\alpha_2}}T_1^{\gamma_1}T_2^{\gamma_2}T_3^{\gamma_3}.
	$$
	On the other hand, the annihilator of $\frac{1}{x_1x_2}$ is given by the left ideal $\T(x_1,x_2)$.
	
	$(ii)$ Follows in a similar way by taking $1$ as the generator.
\end{proof}
\end{proposition}

From this previous proposition we get the isomorphisms of left $\T$-modules
$$
S \cong \T/\T(\partial_1, \partial_2)  \qquad \text{and} \qquad \HL^2(S) \cong \T/\T(x_1, x_2).
$$
\begin{remark}
For any $w \in \HL^2(S)$ we have that $g_i\bullet w = g_i\cdot w$, and so we have that  
$$
\big\{ w \in \HL^2(S) \mid g_1\bullet w=0 \text{ and } g_2\bullet w=0 \big\} = \big\{ w \in \HL^2(S) \mid g_1\cdot w=0 \text{ and } g_2\cdot w=0 \big\}, 
$$
which tells us that we can enlarge $S$ into $\T$ and still recover the same object $\EEQ$ that we are interested in.
\end{remark}

At the moment we have a description of $\EEQ$ as the set of elements in $\HL^2(S)$ vanished by the polynomials $g_1$ and $g_2$, but certainly it would be interesting to have a description as the set of elements in $S$ vanished by certain differential operators.
To achieve this we shall use the Fourier transform (see \cite[Section 5.2]{COUTINHO}).

\begin{definition}
	\label{Fourier_transform}
	By $\Fou$ we denote the automorphism on $\T$ defined by 
	$$
		\Fou(x_i)=\partial_i, \quad \Fou(\partial_i)=-x_i, \quad \Fou(T_i)=T_i. 
	$$
\end{definition}

\begin{notation}
	For the rest of this paper we shall use the notations $L_1 = \Fou(g_1)$ and $L_2 = \Fou(g_2)$.
\end{notation}

\begin{lemma}
	\label{BIGRAD_STRUCT_SOL}
	The $\RR$-vector space $\Sol(L_1,L_2; S)=\big\{ h \in S \mid L_1 \bullet h= 0 \text{ and } L_2 \bullet h= 0 \big\}$ has a structure of $S$-module given by the twisting of the Fourier transform:
	\begin{equation}
		\label{twisting}
			\text{let } f \in S,\; h \in \Sol(L_1,L_2;S) \quad\text{ then we define } \quad f\cdot h = \Fou(f) \bullet h. 		
	\end{equation}
Also it has a bigraded structure induced from S, that is, 
	\begin{equation}
		\label{bigrad_dec_Sol}
		\Sol(L_1, L_2;S) = \bigoplus_{i \ge 0,j \le 0} \Sol(L_1,L_2;S)_{i,j},		
	\end{equation}
	where $\Sol(L_1,L_2;S)_{i,j}=\Sol(L_1,L_2;S) \cap S_{i,-j}$.
	\begin{proof}
		For any $h \in \Sol(L_1,L_2;S)$ we have that $T_ih \in \Sol(L_1,L_2;S)$ and $\partial_i \bullet h \in \Sol(L_1,L_2;S)$, therefore it follows that $\Sol(L_1,L_2;S)$ has a structure of $S$-module given by \autoref{twisting}.
		
		The bigraded decomposition of \autoref{bigrad_dec_Sol} comes from the fact $L_1$ and $L_2$ are bihomogeneous, both with degree $1$ on the $T_i$'s, and degree $\mu$ and $d-\mu$ respectively on the $\partial_i$'s.
		We need to index with non-positive integers $j\le 0$ on the $x$-degree to satisfy the condition $x_i \cdot \Sol(L_1,L_2;S)_{i,j} \subset \Sol(L_1,L_2;S)_{i,j+1}$.
	\end{proof}
\end{lemma}

\begin{notation}
	From now on we shall use the notation $\SE=\TwistSol$ to stress the bigraded $S$-module structure induced on $\Sol(L_1,L_2;S)$ by the twisting of the Fourier transform $\Fou$.
\end{notation}

\begin{theorem}
	\label{isom_Sol}
	Adopt \autoref{data}.
	We have the following isomorphism of bigraded $S$-modules
	$$
	\EEQ \cong \SE(-2,-d+2).
	$$	
	\begin{proof}
		We shall prove that this isomorphism is induced by the Fourier transform $\Fou$.
		We divide the proof into three short steps.
		
		\underline{Step 1.} We define the following two canonical maps
		$$
		\Pi_x : \T \rightarrow \T/\T(\partial_1,\partial_2)\;(\cong S), \quad\qquad\quad \Pi_\partial : \T \rightarrow \T/\T(x_1,x_2) \;(\cong \HL^2(S)).
		$$
		For any $z \in \T$ we have the equivalence
		\begin{equation*}
			\left( 
			\begin{array}{c}
			   z \;\in\; \T(x_1,x_2)\\
			\end{array}			
			\right)
			\Longleftrightarrow
			\left( 
			\begin{array}{c}
				\Fou(z) \;\in\; \T(\partial_1,\partial_2)\\
			\end{array}			
			\right),
		\end{equation*}
		therefore we get an induced isomorphism $\overline{\Fou}:\HL^2(S) \rightarrow S_{\Fou}$ of left $\T$-modules, where $S_{\Fou}$ denotes $S$ twisted by $\Fou$. This isomorphism satisfies
		$$
		\overline{\Fou}(\Pi_\partial(z)) = \Pi_x(\Fou(z)).
		$$
		
		\underline{Step 2.} For any $z \in \T$ we have the following equivalences
		\begin{align*}
				\left( 
				\begin{array}{c}
				g_1\bullet\Pi_\partial(z)=0\\
				g_2\bullet\Pi_\partial(z)=0\\
				\end{array}			
				\right)
				& \Longleftrightarrow
				\left( 
				\begin{array}{c}
				g_1z \;\in\; \T(x_1,x_2)\\
				g_2z \;\in\; \T(x_1,x_2)\\
				\end{array}			
				\right)
				\Longleftrightarrow
				\left( 
				\begin{array}{c}
				\Fou(g_1)\Fou(z) \;\in\; \T(\partial_1,\partial_2)\\
				\Fou(g_2)\Fou(z) \;\in\; \T(\partial_1,\partial_2)\\
				\end{array}			
				\right)
				\Longleftrightarrow\\
				\Longleftrightarrow&
				\left( 
				\begin{array}{c}
				L_1\bullet\Pi_x(\Fou(z))=0\\
				L_2\bullet\Pi_x(\Fou(z))=0\\
				\end{array}			
				\right)
				\Longleftrightarrow
				\left( 
				\begin{array}{c}
				L_1\bullet\overline{\Fou}(\Pi_\partial(z))=0\\
				L_2\bullet\overline{\Fou}(\Pi_\partial(z))=0\\
				\end{array}			
				\right).
		\end{align*}
		Therefore $\overline{\Fou}$ induces an isomorphism of $S$-modules
		$$
		\big\{ w \in \HL^2(S) \mid g_1\bullet w=0 \text{ and } g_2 \bullet w=0\big\}
		\cong 
		\TwistSol.
		$$
		
		\underline{Step 3.} From the definition of $\Fou$ we have that $\overline{\Fou}$ is homogeneous of degree $0$ on the $T_i$'s. On the other hand, we have that $\overline{\Fou}$ makes a shift degree of $2$ in the $x_i$'s since it sends 
		$$
		\frac{1}{x_1^{\alpha_1}x_2^{\alpha_2}}={(-1)}^{\alpha_1+\alpha_2}\frac{\partial_1^{\alpha_1-1}\partial_2^{\alpha_2-1}}{(\alpha_1-1)!(\alpha_2-1)!} \bullet \frac{1}{x_1x_2} \in \HL^2(R)
		$$
		 into 
		 $$
		 \Fou\left({(-1)}^{\alpha_1+\alpha_2}\frac{\partial_1^{\alpha_1-1}\partial_2^{\alpha_2-1}}{(\alpha_1-1)!(\alpha_2-1)!}\right)
		 =
		 \frac{x_1^{\alpha_1-1}x_2^{\alpha_2-1}}{(\alpha_1-1)!(\alpha_2-1)!} \in R.
		 $$
		Then adding the shift degrees $(0,2)$ to \autoref{eqs_whole} we obtain the result.
	\end{proof}
\end{theorem}

\begin{notation}
	Since both $L_1$ and $L_2$ are linear on the $T_i$'s, then we get that $\T(L_1,L_2) \in \LCAT$ and $\T/\T(L_1,L_2) \in \LCAT$.
	We shall denote this last quotient by $Q=\T/\T(L_1, L_2)$.
\end{notation}

Before finishing this section we present an isomorphism of graded $U$-modules that will be the starting point of \hyperref[section5]{Section 5}.

\begin{proposition}
	\label{K_as_Hom}
	Adopt \autoref{data}.
	We have the following isomorphism of graded $U$-modules
	$$
	\EEQ \cong {}^*\Hom_{\T}\big(Q,S\big) (-2).	
	$$
	\begin{proof}
		The following isomorphism of $\RR$-vector spaces
		$$
		\Hom_{\T}\big(\T\big/\T(L_1,L_2),S\big) \cong \big\{ h \in S \mid L_1\bullet h=0 \text{ and } L_2\bullet h=0 \big\}
		= 
		\Sol(L_1,L_2; S)
		$$
		follows in the same way as in \cite[Chapter 6, Theorem 1.2]{COUTINHO}.
		From the discussions of \hyperref[notations_sec]{Subsection 3.1}, we actually have an isomorphism ${}^*\Hom_{\T}(Q, S)\cong \Sol(L_1,L_2;S)$ of graded $U$-modules.
		The shifting of degree follows from \autoref{isom_Sol}.
		\end{proof}
\end{proposition}

\section{The bigraded structure of $\EEQ$ and its relation with $b$-functions}\label{section4}

For organizational purposes we have divided this section into two subsections. 
In the first one, we use the theory of $D$-modules (specifically, the existence of $b$-functions) to bound from above the degree of the polynomial solutions of the system of differential equations $\Sol(L_1,L_2;S)$, then from \autoref{isom_Sol} follows a lower bound in the possible $x$-degree.
In the second subsection, using the local duality theorem for graded modules we prove that this bound it is always strict. 

\subsection{Polynomial solutions}
Our treatment in this subsection will follow \cite[Section 2]{OAKU_GROBNER}, but we need to make some variations since the algorithm given there is restricted to holonomic ideals inside the Weyl algebra.
We shall use \cite{GROBNER_DEFORMATIONS} as our reference regarding Gr\"{o}bner deformations and the algorithmic aspects of $D$-modules.

\begin{notation}
	\label{notations_groebner}
	Adopt \autoref{data}.
	For the rest of the present section we fix the integers $p\ge 2$, $m=\binom{p}{2}$ and $n=\binom{p+1}{2}$. 
	The graded part $\SE_{p-2,*}$
	is given as the solution set of the system of differential equations 
	\begin{equation}
		\label{sys_V}
		V = \big\{h=(h_1,\ldots,h_m) \in R^m \mid [L_1] \bullet h = 0 \text{ and } [L_2] \bullet h = 0  \big\},		
	\end{equation}
	where $[L_i] \in D^{n\times m}$ is an $n\times m$ matrix with entries in $D$ and induced by restricting $L_i$ to the monomials $\mathbf{T}^\gamma$ of degree $\lvert\gamma\rvert=p-2$.
	We join both matrices in a single matrix $H \in D^{2n\times m}$ defined by 
	\begin{equation}
		\label{matrix_H}	
		H = \left( \begin{array}{c}
				L_1\\
				\hline
				L_2
		\end{array} \right),
	\end{equation}
	then equivalently we can write $V = \big\{h=(h_1,\ldots,h_m) \in R^m \mid H \bullet h = 0 \big\}$.
	We define $N \subset D^m$ as the left $D$-module given as the image of $H$, i.e., $N = D^{2n}\cdot H$.
	With $M$ we denote the quotient module $M=D^m/N$, we have an isomorphism $\Hom_D(M, R^m) \cong V$ of $\RR$-vector spaces.
\end{notation}

\begin{example}
	\label{example_sys}
		We give the explicit form of the system of differential equations \autoref{sys_V} in the cases $p=2$ and $p=3$.
		Suppose that $L_1=a_1T_1+a_2T_2+a_3T_3$ and $L_2=b_1T_1+b_2T_2+b_3T_3$.
		For $p=2$ we have that $h=(h_1) \in S_0=R$, and the equations $L_1\bullet h=0$ and $L_2\bullet h=0$ can be expressed as 
		$$
		\left( \begin{array}{c}
			a_1\\
			a_2\\
			a_3
		\end{array} \right)
		\bullet (h_1) = 0
		\quad\text{ and }\quad
		\left( \begin{array}{c}
			b_1\\
			b_2\\
			b_3
		\end{array} \right)
		\bullet (h_1) = 0,
		$$
		and in this case we have that $N$ is actually the left ideal $D(a_1,a_2,a_3,b_1,b_2,b_3)$.
		When $p=3$, we have that $h=(h_1,h_2,h_3)\in S_1=RT_1+RT_2+RT_3$, and sorting the monomials $\mathbf{T}^{\gamma}$ in lexicographical order we get that the equations $L_1\bullet h = 0$ and $L_2 \bullet h=0$ can be expressed as
\end{example}				
		$$
		\begin{blockarray}{cccc}	
			\begin{block}{cccc}
			& T_1 & T_2 & T_3\\	
			\end{block}
			\begin{block}{c(ccc)}
					T_1^2 & a_1 & 0 & 0 \\
					T_1T_2 & a_2 & a_1 & 0 \\			
					T_1T_3 & a_3 & 0 & a_1 \\			
					T_2^2 & 0 & a_2 & 0 \\			
					T_2T_3 & 0 & a_3 & a_2 \\			
					T_3^2 & 0 & 0 & a_3\\
			\end{block}			
		\end{blockarray}
		\bullet
		\left(
		\begin{array}{c}
			h_1\\
			h_2\\
			h_3
		\end{array}
		\right)
		= 0
		\quad\text{ and }\quad
		\begin{blockarray}{cccc}	
		\begin{block}{cccc}
		& T_1 & T_2 & T_3\\	
		\end{block}
		\begin{block}{c(ccc)}
		T_1^2 & b_1 & 0 & 0 \\
		T_1T_2 & b_2 & b_1 & 0 \\			
		T_1T_3 & b_3 & 0 & b_1 \\			
		T_2^2 & 0 & b_2 & 0 \\			
		T_2T_3 & 0 & b_3 & b_2 \\			
		T_3^2 & 0 & 0 & b_3\\
		\end{block}			
		\end{blockarray}
		\bullet
		\left(
		\begin{array}{c}
		h_1\\
		h_2\\
		h_3
		\end{array}
		\right)
		= 0.
		$$	
%\end{example}

\begin{proposition}
	\label{M_is_holonomic}
	The left $D$-module $M=D^m/N$ is holonomic.
	\begin{proof}
			From the exact sequence \autoref{Koszul_comp} of \autoref{exactness_lemma} we take the graded part $p$ in the $T_i$'s, which gives the following exact sequence of left $D$-modules
			\begin{equation}
				\label{graded_part_Koszul}
					0 \rightarrow \T_{p-2} \xrightarrow{\cdot A} \T_{p-1}^2 \xrightarrow{\cdot B} \T_p \rightarrow Q_p \rightarrow 0,
			\end{equation}
			where $A = \left(-{[L_2]}^T \mid {[L_1]}^T \right)$ and $[L_i]^T$ represent the transpose of the matrix $[L_i]$ from \autoref{notations_groebner} \footnote{With the same notations of \autoref{example_sys}, when $p=2$ we have $A=\left(\begin{array}{lll|ccc}
			 -b_1& -b_2 & -b_3 & a_1 & a_2 & a_3
			\end{array}
			\right) \in D^{1 \times 6}$.}.
							
			Applying the functor $\Hom_D(-,D)$ we obtain the following complex of right $D$-modules 
			$$
			0 \rightarrow {(\T_p)}^T \xrightarrow{B\cdot} {(\T_{p-1}^2)}^T \xrightarrow{A\cdot} {(\T_{p-2})}^T \rightarrow 0,
			$$
			where the cokernel of the last map $A\cdot$  is  $\Ext_D^2(Q_p, D)$, but from the form of $A$ this coincides with the standard transposition $\tau(M)$ of $M$, that is, $\Ext_D^2(Q_p,D)\cong \tau(M)$.
			Finally, from \cite[Lemma 7.3, page 73]{BJORK}	we have that $\Ext_D^2(Q_p,D)$ is a holonomic right $D$-module, and this clearly implies the holonomicity of $\tau(M)$ and $M$.					
	\end{proof}
\end{proposition}

We recall the notion of Gr\"obner deformations (see \cite[Section 1.1]{GROBNER_DEFORMATIONS}). For a given generic weight $w=(w_1,w_2)$ and an element $\ell=\sum_{\alpha,\beta} c_{\alpha,\beta} x^{\alpha}\partial^{\beta} \in D$, we denote by $\iniTerm_{(-w,w)}(\ell)$ the initial form of $\ell$ with respect to $w$ and it is defined as
\begin{equation}
	\label{initial_Form}
	\iniTerm_{(-w,w)}(\ell) = \sum_{-\alpha\cdot
	 w + \beta\cdot w \text{ is maximum}} c_{\alpha,\beta}x^{\alpha}\partial^{\beta}.
\end{equation}
We are interested in the weight $w=(-1,-1)$ that makes $\deg(x_i)=1$ and $\deg(\partial_i)=-1$, thus we shall drop the subscript $(-w,w)$ and the definition of the initial form \autoref{initial_Form} turns into
\begin{equation}
\label{my_initial_Form}
\iniTerm(\ell) = \sum_{\lvert\alpha
	- \beta\rvert \text{ is maximum}} c_{\alpha,\beta}x^{\alpha}\partial^{\beta}.
\end{equation}

\begin{definition}\textit{\cite[Corollary 1.1.2, Definition 1.1.3]{GROBNER_DEFORMATIONS}}
	Let $J \subset D$ be a left ideal, then the $\RR$-vector space
	$$
	\iniTerm(J)=\RR\cdot\Big\{\iniTerm(\ell) \mid \ell \in J \Big\}
	$$
	is a left ideal in $D$ and it is called the initial ideal of $J$.
\end{definition}

\begin{definition}\textit{\cite[Definition 5.1.1]{GROBNER_DEFORMATIONS}}
	Let $J \subset D$ be a holonomic left ideal. The  elimination ideal 
	\begin{equation}
		\label{elim_ideal}
		\iniTerm(J) \;\cap\; \RR[-x_1\partial_1 - x_2\partial_2]
	\end{equation}
	is principal in the univariate polynomial ring $\RR[s]$, where $s=-x_1\partial_1-x_2\partial_2$.  
	The generator $b_J(s)$ of the principal ideal \autoref{elim_ideal} is called the $b$-function of $J$.
\end{definition}

An important fact is that the $b$-function of a holonomic ideal is a non-zero polynomial (see e.g. \cite[Theorem 5.1.2]{GROBNER_DEFORMATIONS}).
Now we present suitable definition for the $b$-function of a left $D$-module, which is essentially the same as the one given in \cite[Section 4]{OAKU}
(see \cite[Lemma 4.2]{OAKU}).

\begin{definition}
	\label{bFunction_module}
	Let $M^{'}$ be a holonomic left $D$-module given as the quotient module $M^{'}=D^r/N^{'}$.
	For each $i = 1,\ldots,r$ with the canonical projection $\pi_i:D^r\rightarrow D$ of $D^r$ onto the $i$-th component $e_i$, we define a left $D$-ideal 
	$$
	J_i=\pi_i(N^{'} \;\cap \;D\cdot e_i)=\big\{\ell \in D \mid (0,\ldots,\underbrace{\ell}_{i\text{-th}},\ldots,0) \in N^{'}  \big\}.
	$$
	Then the $b$-function of $M^{'}$ is given as the least common multiple of the $b$-functions of the $D$-ideals $J_i$, that is,
	$$
	b_{M^{'}}(s) ={\normalfont\text{LCM}}_{i=1,\ldots,r} \;\big(b_{J_i}(s)\big).
	$$ 
\end{definition}

In this previous definition for each $i=1,\ldots,r$, the canonical injection $D/J_i \hookrightarrow D^r/N^{'}$ implies that each ideal $J_i$ is holonomic, and so we get that the $b$-function of a holonomic module is a non-zero polynomial.
Before proving the main result of this subsection, we recall an easy but important lemma.

\begin{lemma}
	\label{Euler_formula}
	Let $P(s) \in\RR[s]$ be a polynomial in $s=-x_1\partial_1-x_2\partial_2$ and let $f \in R$ be a homogeneous polynomial of degree $\deg(f)=k$, then we have $P(s) \bullet f = P(-k)f$.
	\begin{proof}
		It follows from Euler's formula $(x_1\partial_1+x_2\partial_2)\bullet f = kf$.
	\end{proof}
\end{lemma}

\begin{theorem}
	Adopt \autoref{data}.
	Consider the $b$-function $b_M(s)$ of the holonomic $D$-module $M$ defined in \autoref{notations_groebner}.
	For any integer $q$, if $b_M(-d+2+q)\neq 0$ then we have that $\EEQ_{p,q}=0$.
	\begin{proof}
		Suppose by contradiction  that $\EEQ_{p,q}\neq 0$, then from \autoref{isom_Sol} there exists $0 \neq h \in \SE_{p-2,-k}$ where $-k=-d+2+q$.
		Indexing this element as in \autoref{notations_groebner} we have a non-zero polynomial vector $h=(h_1,\ldots,h_m) \in V$ where each polynomial $h_i$ has degree $\deg(h_i)=k$.
		
		For each each $i=1,\ldots,m$, let $b_{J_i}(s)$ be the $b$-function corresponding to the left $D$-ideal $J_i=\pi_i(N \;\cap \;D\cdot e_i)$.
		Then we have $b_{J_i}(s)\cdot e_i \bullet h = 0$, which implies $b_{J_i}(s) \bullet h_i=0$. 
		Using \autoref{Euler_formula} we get $b_{J_i}(-k)h_i=0$, but since $b_{J_i}(-k) \neq 0$ then we have $h_i=0$.
		Finally, we have obtained the contradiction $h=0$.
	\end{proof}
\end{theorem}

\begin{corollary}
	\label{divide_b_function}
		Adopt \autoref{data}.
		Let $q$ be the lowest possible $x$-degree for an element in the graded part $\EEQ_{p,*}$, that is, $\EEQ_{p,q} \neq 0$ and $\EEQ_{p,q-1}=0$.
		Then the polynomial $s(s+1)\cdots(s+d-2-q)$ divides the $b$-function $b_M(s)$.
		\begin{proof}
			Follows from the contrapositive of the previous theorem.
		\end{proof}
\end{corollary}

\subsection{The equality}

In this subsection we shall prove that the approximation given above it is actually strict.

\begin{lemma}
	\label{belong_lemma}
	For any $k \ge 0$ we have the identity 
	$$
	s(s+1)\cdots(s+k)={(-1)}^{k+1}\sum_{j=0}^{k+1} \binom{k+1}{j} x_1^jx_2^{k+1-j}\partial_1^j\partial_2^{k+1-j}.
	$$
	Thus, we have that 
	\begin{enumerate}[(i)]
		\item $s(s+1)\cdots(s+k) \in D(\partial_1, \partial_2)^{k+1}$, where $D(\partial_1,\partial_2)^{k+1}$ denotes the left $D$-ideal generated by the elements $\{\partial_1^{\beta_1}\partial_2^{\beta_2}\mid \beta_1+\beta_2 = k+1 \}$;
		\item $s(s+1)\cdots(s+k)$ is homogeneous, that is $$
		\iniTerm\big(s(s+1)\cdots(s+k)\big)=s(s+1)\cdots(s+k).$$
	\end{enumerate}
	\begin{proof}
		We proceed by induction on $k$. 
		For $k=0$ it is clear since $s=-x_1\partial_1-x_2\partial_2$. 
		
		First we prove the identity
		 $x_i^{\beta_i}\partial_i^{\beta_i}(x_i\partial_i - \beta_i)=x_i^{\beta_i+1}\partial_i^{\beta_i+1}$  using induction on $\beta_i$.
		For $\beta_i=0$ it is vacuous, thus we assume that $\beta_i > 0$ and that the statement holds for any non-negative integer smaller than $\beta_i$.
		Hence, we have the equalities 
		\begin{align*}
			x_i^{\beta_i}\partial_i^{\beta_i}(x_i\partial_i-\beta_i)&=x_i^{\beta_i}\partial_i^{\beta_i-1}(\partial_ix_i\partial_i-\beta_i\partial_i)\\
			&=x_ix_i^{\beta_i-1}\partial_i^{\beta_i-1}(x_i\partial_i-(\beta_i-1))\partial_i=x_ix_i^{\beta_i}\partial_i^{\beta_i}\partial_i=x_i^{\beta_i+1}\partial_i^{\beta_i+1}.
		\end{align*}
		Then we can obtain that
		\begin{align*}
			x_1^jx_2^{k+1-j}\partial_1^j\partial_2^{k+1-j}(s+k+1)=(-1)\big(x_1^{j+1}x_2^{k+1-j}\partial_1^{j+1}\partial_2^{k+1-j} + x_1^jx_2^{k+2-j}\partial_1^j\partial_2^{k+2-j} \big),
		\end{align*}
		and using $\binom{k+1}{j}+\binom{k+1}{j-1}=\binom{k+2}{j}$ the proof of the lemma follows similarly to the usual binomial theorem.
\end{proof}
\end{lemma}
From \autoref{notations_groebner} we can define the matrix $F = \Fou(H) = (\Fou(H_{i,j})) \in R^{2n\times m}$ (where $m=\binom{p}{2}, n=\binom{p+1}{2}$), that is, the $2n\times m$ matrix with entries in $R$ obtained after applying the Fourier  transform to each entry of the matrix $H$.
In a similar way to \autoref{notations_groebner}, we define the graded $R$-module $L = R^{m}/(R^{2n}\cdot F)$ (all the rows of $F$ are homogeneous of degree $\mu$ or degree $d-\mu$).

Since $\{g_1,g_2\}$ is a regular sequence in $S$ (see the proof of \autoref{exactness_lemma}), by restricting the Koszul complex $K_{\bullet}(g_1,g_2)$ to the graded part $p$, the module $\Sym_{p}(I)$ gets the graded free resolution
$$
0 \rightarrow {R(-d)}^{\binom{p}{2}} \rightarrow
{R(-\mu)}^{\binom{p+1}{2}}\oplus{R(-d+\mu)}^{\binom{p+1}{2}}\rightarrow R^{\binom{p+2}{2}} \rightarrow \Sym_{p}(I) \rightarrow 0. 
$$
 Similarly to \autoref{M_is_holonomic}, when we apply $\Hom_R(-, R)$ we get a complex
\begin{equation}
	\label{res_L}
 0 \rightarrow {R}^{\binom{p+2}{2}} \rightarrow {R(\mu)}^{\binom{p+1}{2}} \oplus {R(d-\mu)}^{\binom{p+1}{2}} \rightarrow {R(d)}^{\binom{p}{2}} \rightarrow 0, 
\end{equation}
where the cokernel of the map on the right is the graded $R$-module ${}^*\Ext_R^2(\Sym_{p}(I), R)$. Making a shift degree of $-d$ on the modules of \autoref{res_L} gives us a complex that has the module $L$ as the cokernel of the map on the right.
Therefore, we have an isomorphism $L(d)\cong {}^*\Ext_R(\Sym_{p}(I), R)$ of graded $R$-modules.

Now, as an application of the local duality theorem in the graded case (see e.g. \cite[Section 14.4]{Brodmann_Sharp_local_cohom} or \cite[Section 3.6]{BRUNS_HERZ}) we can prove our sought equality.
\begin{theorem}
	\label{expression_bFunction}
	Adopt \autoref{data}.
	Let $b_M(s)$ be the $b$-function of the holonomic module $M$ defined in \autoref{notations_groebner} and let $q$ be the lowest possible $x$-degree for an element in the graded part $\EEQ_{p,*}$.
	Then 
	$$
	b_M(s)=s(s+1)\cdots(s+d-2-q).
	$$
	\begin{proof}
		From \autoref{divide_b_function} we already know that $s(s+1)\cdots(s+d-2-q) \mid b_M(s)$, then it will be enough for us to prove that for each $i=1,\ldots,m$ we have 
		$$
		s(s+1)\cdots(s+d-2-q) \in \iniTerm(J_i)\;\cap\;\RR[s],
		$$	
		where $J_i = \pi_i(N \cap D\cdot e_i)$.		
		
		Let $a=\text{end}(L)=\max\{k \mid L_k \neq 0 \}$ (since $L$ is a finite length module), then for any $x_1^{\alpha_1}x_2^{\alpha_2}$ with $\alpha_1+\alpha_2=a+1$ we have that 
		$$
		x_1^{\alpha_1}x_2^{\alpha_2}e_i=(0,\ldots,\underbrace{x_1^{\alpha_1}x_2^{\alpha_2}}_{i\text{-th}},\ldots,0) \in R^{2n}\cdot F,
		$$
		where $i=1,\ldots,m$ and by an abuse of notation $e_i$ also represents the $i$-th component of the free $R$-module $R^m$.
		Applying the Fourier transform and using \autoref{belong_lemma} we obtain that 
		$$
		s(s+1)\cdots(s+a) \in \iniTerm(J_i) \;\cap\;\RR[s]
		$$
		for each $i=1,\ldots,m$.
		From the local duality theorem in the graded case, we get the following isomorphisms of graded $R$-modules
		$$
		\EEQ_{p,*}=\HL^0(\Sym_{p}(I)) \cong {}^*\Hom_{\RR}\big({}^*\Ext_R^2(\Sym_{p}(I), R(-2)), \RR \big)
		\cong
		{}^*\Hom_{\RR}\big(L(d-2), \RR \big).
		$$
		Since the grading of ${}^*\Hom_{\RR}\big(L(d-2), \RR \big)$ is given by $${{}^*\Hom_{\RR}\big(L(d-2), \RR \big)}_i={}^*\Hom_{\RR}\big({L(d-2)}_{-i}, \RR \big),$$ we have that $a=d-2-q$, and so the statement of theorem follows.
	\end{proof}
\end{theorem}

\section{Computing Hom with duality}\label{section5}

The aim of this section is to compute ${}^*\Hom_{\T}\big(Q, S\big)$ (where $Q=\T/\T(L_1,L_2)$) by means of some duality that was previously exploited in \cite{HOM_DUALITY}.
In the general Weyl algebra $A_n(\RR)$, for two holonomic left $A_n(\RR)$-modules $M$ and $N$ we have the following duality (see e.g.  \cite[Proposition 4.14, page 58]{BJORK} or  \cite[Theorem 2.1]{HOM_DUALITY})
$$
\Ext_{A_n(\RR)}^i(M, N) \cong \Tor_{n-i}^{A_n(\RR)}\big(\Ext_{A_n(\RR)}^n(M,A_n(\RR)), \;N \big),
$$
which is one of the main tools used in \cite{HOM_DUALITY}.
Unfortunately we want to work over our previously defined algebra $\T$ and for this we will have to make some variations. Nevertheless, we can achieve the following duality in our case.

\begin{theorem}
	\label{duality}
	For any $i$ we have the following isomorphism of graded $U$-modules (see \hyperref[notations_sec]{Subsection 3.1})
	$$
	{}^*\Ext_{\T}^i(Q, S) \cong {}^*\Tor_{2-i}^{\T}\big({}^*\Ext_{\T}^2(Q, \T), \;S \big).
	$$
\end{theorem}

To prove this duality we shall use \cite[Chapter 2]{BJORK} as our main source. We start by defining a Bernstein filtration on $\T$ and exploiting the induced graded ring.

\begin{defprop}
	For any $i \ge 0$ we define the $\RR$-vector space $F_i$ which is generated by the set of monomials $\{ x_1^{\alpha_1}x_2^{\alpha_2}\partial_1^{\beta_1}\partial_2^{\beta_2}T_1^{\gamma_1}T_2^{\gamma_2}T_3^{\gamma_3} \mid\; \mid\mathbf{\alpha}\mid + \mid\mathbf{\beta}\mid + \mid\mathbf{\gamma}\mid \;\le\; i \}$, and we denote $F_{-1}=0$.
	Since we have 
	\begin{enumerate}[(1)]
		\item $0 = F_{-1} \subset F_0 \subset F_1 \subset F_2 \subset \cdots \subset \T$,
		\item $\T = \bigcup_{i\ge 0} F_i$,
		\item $F_i \cdot F_j \subset F_{i+j}$,
	\end{enumerate}
	then $F=\{F_i\}$ is a filtration of $\T$. 
	With this filtration we induce the associated graded ring
	{\normalfont $\gr(\T)=\bigoplus_{i\ge 0} F_i/F_{i-1}$}, which is isomorphic to a polynomial ring of $7$ variables with coefficients in $\RR$. We shall use the notation {\normalfont $$\gr(\T) \cong T := \RR[x_1,x_2,\delta_1,\delta_2, T_1,T_2,T_3],$$} where we get a canonical map $\sigma: \T \rightarrow T$ given by $\sigma(x_i)=x_i$, $\sigma(\partial_i) = \delta_i$ and $\sigma(T_i)=T_i$.
	\begin{proof}
		See \cite[Proposition 2.2, page 4]{BJORK} or \cite[Theorem 3.1, page 58]{COUTINHO}.
	\end{proof}
\end{defprop}

We shall denote by $q_1=\sigma(L_1)$ and $q_2=\sigma(L_2)$ the elements in $T$ corresponding to $L_1$ and $L_2$. Here we have, that $q_1$ and $q_2$ are bihomogeneous polynomials which are linear on the $T_i$'s, and have degree $\mu$ and $d-\mu$ on the $\delta_i$'s respectively.
But from the graded structure of $T$, we only see them as homogeneous polynomials having degree $\mu+1$ and $d-\mu+1$ respectively.

A filtration on a left $\T$-module $M$ consists of an increasing sequence of finite dimensional $\RR$-subspaces $0=\Gamma_{-1}\subset \Gamma_0 \subset \Gamma_1 \subset \Gamma_2 \subset \cdots$ satisfying $\bigcup \Gamma_i=M$ and the inclusions $F_i\cdot\Gamma_j \subset \Gamma_{i+j}$ for all $i$ and $j$. 
With a filtration we get the associated graded $T(=\gr(\T))$-module $\gr_\Gamma(M)=\bigoplus_{i\ge 0} \Gamma_i/\Gamma_{i-1}$. 
We say that $\Gamma=\{\Gamma_i\}$ is a good filtration if $\gr_\Gamma(M)$ is a finitely generated $T$-module.
Using a good filtration we can define a Hilbert-Samuel function, and so we can get a notion of dimension for left $\T$-modules.

\begin{defprop}
	Given a good filtration $\Gamma = \{\Gamma_i\}$ for a finitely generated left $\T$-module $M$, then there exists a polynomial $\chi_M^\Gamma(t)=a_dt^d+ \cdots + a_1t+a_0$ with rational coefficients such that $\dim_{\RR}(\Gamma_t) = \chi_M^\Gamma(t)$ when $t \gg 0$. 
	The integer $d$ is independent of the good filtration chosen, and we define $d(M)=d$ as the Bernstein dimension of $M$.
	\begin{proof}
		See \cite[Section 1.3]{BJORK} or \cite[Chapter 9]{COUTINHO}.
	\end{proof}
\end{defprop}

Since the Hilbert-Samuel function of $T$ is given by $\binom{t+7}{7}$, thus we have $d(\T)=7$.
Now we want to study the left $\T$-module $Q=\T/\T(L_1,L_2) \in \LCAT$, and we begin by proving that the Koszul complex gives a free resolution for it in $\LCAT$. 

\begin{proposition}
	\label{exactness_lemma}
	The following statements hold.
	\begin{enumerate}[(i)]
		\item The dimension of $Q$ is $d(Q)=5$.
		\item The following  Koszul complex in $\LCAT$ is exact
		\begin{equation}
			\label{Koszul_comp}
		 \mathcal{A}: \quad	0 \rightarrow \T(-2) \xrightarrow{\cdot\left[\begin{smallmatrix}
			-L_2,L_1
			\end{smallmatrix}\right]} \T(-1)^2 \xrightarrow{\cdot\left[\begin{smallmatrix}
			L_1\\
			L_2
			\end{smallmatrix}\right]}
		\T \rightarrow Q \rightarrow 0.
		\end{equation}
		
	\end{enumerate}
	\begin{proof}
		$(i)$ 
		The module $Q$ being a quotient of $\T$ automatically gets a natural good filtration given by $F_i/(F_i\cap\T(L_1,L_2))$ and then we get $\gr(Q)=\bigoplus_{i\ge0}F_i/(F_{i-1}+F_i\cap\T(L_1,L_2))\cong T/(q_1,q_2)$.
		
		Let $B=\RR[x_1,x_2,\delta_1,\delta_2]$, and let $h_i$ be the polynomials in $B$ obtained from $f_i$ by making the substitution $x_i\mapsto \delta_i$, i.e., $h_i=\sigma(\Fou(f_i))$. 
		By the Hilbert-Burch theorem we have that $J=(h_1,h_2,h_3)\subset B$ is a perfect ideal of height two, and making the substitution $x_i \mapsto \delta_i$ in the resolution \autoref{HilbertBurchRes} of $I$ we get a resolution 
		$$
		0 \rightarrow B^2 \xrightarrow{\varPhi} B^3 \xrightarrow{\left[\begin{smallmatrix}
			h_1,h_2,h_3
			\end{smallmatrix}\right]} J \rightarrow 0,
		$$
		where $\left[q_1,q_2\right]=\left[T_1,T_2,T_3\right]\cdot\varPhi$.
		Hence $T/(q_1,q_2)=\Sym(J)$, and from \cite[Corollary 2.2]{SIMIS_VASC_SYM_INT} the Krull dimension is given by $\dim(T/(q_1,q_2))=\dim(\Sym(J))=\dim(B)+\rank(J)=4+1=5$. Finally, this coincides with the degree of the Hilbert-Samuel function, i.e., $d(Q)=\dim(T/(q_1,q_2))=5$ (see e.g. \cite[Theorem 13.4]{MATSUMURA}).

		$(ii)$ The shifting of degrees in \autoref{Koszul_comp} are clear since $L_1$ and $L_2$ are both linear on the $T_i$'s, then it will be enough to prove exactness of \autoref{Koszul_comp} just in the category $\T$ (i.e., forgetting the graded structures induced in \hyperref[notations_sec]{Subsection 3.1}).
		So, inside this proof, to avoid confusions the only additional structure that we assume on $\T$ is the Bernstein filtration and the induced graded ring $T$.

		Since $\T$ is non-commutative we should check that \autoref{Koszul_comp} is even a complex, but fortunately $L_1$ and $L_2$ are only defined in the $\partial_i$ and $T_i$ variables, and so $L_1L_2-L_2L_1=0$.
		
		The complex \autoref{Koszul_comp} induces the following graded Koszul complex in $T$
		\begin{equation*}			
			0 \rightarrow T(-d-2) \xrightarrow{\left[\begin{smallmatrix}
				-q_2,q_1
			\end{smallmatrix}\right]} T(-\mu-1)\oplus T(-d+\mu-1) \xrightarrow{\left[\begin{smallmatrix}
				q_1\\
				q_2
			\end{smallmatrix}\right]}
			T \rightarrow T/(q_1,q_2) \rightarrow 0.			
	\end{equation*}
	Using that $\dim(T/(q_1,q_2))=5$ we get that $(q_1,q_2)$ is a $T$-regular sequence (see e.g. \cite[Theorem 17.4]{MATSUMURA}) and so this new complex is exact. 
	Finally, \cite[Lemma 3.13, page 46]{BJORK} implies that \autoref{Koszul_comp} is exact.		
	\end{proof}
\end{proposition}

\begin{corollary}
	\label{vanishing_Ext}
	For any $j\neq 2$ we have {\normalfont ${}^*\Ext_{\T}^j(Q,\T)=0$}, and  ${}^*\Ext_{\T}^2(Q,\T)\neq 0$.
	\begin{proof}
		Since \autoref{Koszul_comp} is a free resolution of $Q$ we clearly have ${}^*\Ext_{\T}^j(Q,\T)=0$ for $j>2$.
		On the other hand from \cite[Theorem 7.1, page 73]{BJORK} we have that $j(Q) + d(Q) = 7$, where $j(Q)=\inf\{k \mid {}^*\Ext_{\T}^k(Q,\T)\neq 0\}$.
		Since $d(Q)=5$, then $j(Q)=2$ and the statement of the corollary follows.
	\end{proof}
\end{corollary}

Now we are ready to prove the duality that we claimed in the beginning of this section.

\begin{proof}[Proof of \autoref{duality}]
	A resolution of $S$ in $\LCAT$ is given by the Koszul complex 
	
	\begin{equation}
	\label{sec_Koszul_complex}
	\mathcal{B}:\quad 0 \rightarrow \T \xrightarrow{\cdot\left[\begin{smallmatrix}
		-\partial_2,\partial_1
		\end{smallmatrix}\right]} \T^2 \xrightarrow{\cdot\left[\begin{smallmatrix}
		\partial_1\\
		\partial_2
		\end{smallmatrix}\right]}
	\T \rightarrow S \rightarrow 0.		
	\end{equation}
	
	We define the following third quadrant double complex ${}^*\Hom_{\T}(\AAA, \T) \otimes_\T \BBB$,
	
	\begin{center}	
		\begin{tikzpicture}
		\matrix (m) [matrix of math nodes,row sep=1.5em,column sep=2em,minimum width=2em, text height=1.5ex, text depth=0.25ex]
		{
			{}^*\Hom_{\T}(\AAA_2,\T)\otimes_\T \BBB_2 & {}^*\Hom_{\T}(\AAA_1,\T)\otimes_\T \BBB_2 & {}^*\Hom_{\T}(\AAA_0, \T)\otimes_\T \BBB_2 \\
			{}^*\Hom_{\T}(\AAA_2,\T)\otimes_\T \BBB_1 & {}^*\Hom_{\T}(\AAA_1,\T)\otimes_\T \BBB_1 & {}^*\Hom_{\T}(\AAA_0, \T)\otimes_\T \BBB_1 \\
			{}^*\Hom_{\T}(\AAA_2,\T)\otimes_\T \BBB_0 & {}^*\Hom_{\T}(\AAA_1,\T)\otimes_\T \BBB_0 & {}^*\Hom_{\T}(\AAA_0, \T)\otimes_\T \BBB_0. \\
		};
		\path[-stealth]
		(m-1-2) edge (m-1-1)
		(m-1-3) edge (m-1-2)
		(m-2-2) edge (m-2-1)
		(m-2-3) edge (m-2-2)
		(m-3-2) edge (m-3-1)
		(m-3-3) edge (m-3-2)
		(m-1-1) edge (m-2-1)
		(m-2-1) edge (m-3-1)
		(m-1-2) edge (m-2-2)
		(m-2-2) edge (m-3-2)
		(m-1-3) edge (m-2-3)
		(m-2-3) edge (m-3-3)
		;				
		\end{tikzpicture}	
	\end{center}

	Thanks to our construction of $\LCAT$ and $\RCAT$, we have that this double complex fits naturally in the category of graded $U$-modules, that is, all its elements are graded $U$-modules and all its maps are homogeneous homomorphisms of graded $U$-modules.
	
	Since each ${}^*\Hom_{\T}(\AAA_j, \T) \in \RCAT$ is a free module then by computing homology on each column we get that the only row that does not vanish is the last one. 
	On the other hand, from \autoref{vanishing_Ext} we have that when we compute homology on each row the only column that does not vanish is the leftmost one.
	
	Therefore the spectral sequence determined by the first filtration is given by
	\begin{equation*}
		{}^{\text{I}}E_2^{p,q}=\begin{cases}
		{}^*\Ext_{\T}^p\big(Q,S\big)  \quad \text{if } q = 2,\\
		0 \qquad \qquad\quad\text{otherwise,}
		\end{cases}
	\end{equation*}
	and the spectral sequence determined by the second filtration is given by
	\begin{equation*}
		{}^{\text{II}}E_2^{p,q}=\begin{cases}
		{}^*\Tor_{2-q}^\T({}^*\Ext_{\T}^2\big(Q,\T),S\big) \quad \text{if } p = 2,\\
		0 \qquad\qquad\qquad\qquad \qquad\text{otherwise.}
		\end{cases}
	\end{equation*}

	From the fact that both spectral sequences collapse we get the following isomorphisms of graded $U$-modules
	\begin{equation*}
		{}^{\text{I}}E_2^{i,2} \cong H^{i+2}(\text{Tot}({}^*\Hom_{\T}\big(\AAA, \T) \otimes_\T \BBB)\big) \cong {}^{\text{II}}E_2^{2,i},
	\end{equation*}	and so we obtain the duality of the theorem.
\end{proof}

\begin{theorem}
	\label{deRham_cohomology}
	Adopt \autoref{data}.
	Then we have the following isomorphism of graded $U$-modules
	$$\EEQ \cong \Hdr^0(Q) = \{w \in Q \mid \partial_1 \bullet w= 0 \text{ and } \partial_2 \bullet w= 0 \}.
	$$
	In particular, for any integer $p$ we have an isomorphism of $\RR$-vector spaces
	$$
	\EEQ_{p,*} \cong \Hdr^0(Q_p) = \{w \in Q_p \mid \partial_1 \bullet w= 0 \text{ and } \partial_2 \bullet w= 0 \}.
	$$
	\begin{proof}
		From the resolution \autoref{Koszul_comp} of $Q$ we get the following complex in $\RCAT$
		$$
		{}^*\Hom_\T(\AAA,\T): \quad	0 \rightarrow \T \xrightarrow{\left[\begin{smallmatrix}
			L_1\\
			L_2
			\end{smallmatrix}\right]\cdot} \T(1)^2 \xrightarrow{\left[\begin{smallmatrix}
			-L_2, L_1
			\end{smallmatrix}\right]\cdot}
		\T(2)  \rightarrow 0,			
		$$
		then computing the second cohomology of this complex gives that ${}^*\Ext_{\T}^2(Q,\T) \cong \left(\T/(L_1,L_2)\T\right)(2)$, where  $\T/(L_1,L_2)\T=\tau(Q)$ is the standard transposition of $Q$.		
		
		Since the Koszul complex \autoref{sec_Koszul_complex} gives a resolution of $S$, then computing the second homology of the Koszul complex $\tau(Q)(2) \otimes_\T \BBB$  gives the following isomorphisms of graded $U$-modules
		\begin{align*}
		{}^*\Tor_2^\T({}^*\Ext_\T^2(Q,\T), S) &\cong H_2\big( \tau(Q)(2) \otimes_\T \BBB \big)\\
		&\cong
		\{w \in \tau(Q)(2) \mid  w\bullet \partial_1= 0 \text{ and }  w\bullet \partial_2= 0 \}.
		\end{align*}
		From the fact that $\tau(T_i)=T_i$, we have an isomorphism of graded $U$-modules
		$$		
		\{w \in \tau(Q)(2) \mid  w\bullet \partial_1= 0 \text{ and }  w\bullet \partial_2= 0 \}
		\cong		
		\{w \in Q(2) \mid \partial_1 \bullet w = 0 \text{ and }  \partial_2\bullet w = 0 \},
		$$
		then from \autoref{K_as_Hom} and \autoref{duality} we get the following isomorphisms of graded $U$-modules
		\begin{align*}
			\EEQ&\cong {}^*\Hom_\T(Q, S)(-2)\\
			&\cong{}^*\Tor_2^\T({}^*\Ext_\T^2(Q,\T), S)(-2)\\
			&\cong\{w \in Q \mid \partial_1 \bullet w = 0 \text{ and }  \partial_2\bullet w = 0 \},
		\end{align*}
		that imply the statement of the theorem.
	\end{proof} 
\end{theorem}

\section{Examples and computations}\label{section6}	

	In this short section we show a simple script in \textit{Macaulay2} \cite{M2} that we have implemented to compute the $b$-function of each $D$-module $M$ from \autoref{notations_groebner}.
	Actually, we have to say that an enormous number of examples and computations led us to believe the equality of \autoref{expression_bFunction} in the first place.	
	\vspace*{.4cm}

	\begin{spacing}{0.73}
		{   \fontfamily{cmtt}
			\selectfont
			\hspace*{-.43cm}\verb|needsPackage "Dmodules"|\\					
			\verb|bFunctionRees = (I, p) -> (|\\
			\hspace*{.5cm} \verb|R := ring I;| \\
			\hspace*{.5cm} \verb|W := makeWeylAlgebra R;|\\
			\hspace*{.5cm} \verb|T := W[T1, T2, T3], U := QQ[Z1, Z2, Z3];|\\
			\hspace*{.5cm} \verb|A := Fourier (map(W, R, {(vars W)_(0,0),(vars W)_(0,1)})) (res I).dd_2;|\\
			\hspace*{.5cm} \verb|L := matrix{{T1, T2, T3}} * A;|\\ 
			\hspace*{.5cm} \verb|L1 := L_(0, 0), L2 := L_(0, 1);|\\   
			\hspace*{.5cm} \verb|src := flatten entries (map(T, U, {T1, T2, T3})) basis(p - 2, U);|\\
			\hspace*{.5cm} \verb|dest := flatten entries (map(T, U, {T1, T2, T3})) basis(p - 1, U);|\\ 
			\hspace*{.5cm} \verb|m := #src, n := #dest;|\\
			\hspace*{.5cm} \verb|H := mutableMatrix(W, m, 2 * n);|\\	        
			\hspace*{.5cm} \verb|for i from 0 to m - 1 do (|\\   
			\hspace*{1cm} \verb|mult1 := src#i * L1;|\\
			\hspace*{1cm} \verb|mult2 := src#i * L2;|\\
			\hspace*{1cm} \verb|for j from 0 to n - 1 do (|\\	 			
			\hspace*{1.5cm} \verb|R1 := mult1 // gens ideal(dest#j);|\\
			\hspace*{1.5cm} \verb|R2 := mult2 // gens ideal(dest#j);|\\
			\hspace*{1.5cm} \verb|H_(i, j) = (map(W, T, {1, 1, 1})) R1_(0, 0);|\\
			\hspace*{1.5cm} \verb|H_(i, j + n) = (map(W, T, {1, 1, 1})) R2_(0, 0);|\\   
			\hspace*{1cm} \verb|);|\\
			\hspace*{.5cm} \verb|);|\\
			\hspace*{.5cm} \verb|bM := bFunction(coker matrix H, {-1,-1}, toList(m:0));|\\	
			\hspace*{.5cm} use R;	\\
			\hspace*{.5cm} bM     \\\
			)\\
		}
	\end{spacing}
	
We will carry out a couple of examples to show how we can use \autoref{expression_bFunction} to deduce the bigraded structure of $\EEQ$.	
We can save the previous code in a file that we will call \textrm{``bFunctionRees.m2''}

\begin{example}
	Let $I=(x^{5}, x^2y^3, y^{5}) \subset \mathbb{Q}[x,y]$, then from \cite{MONOMIAL_CARLOS}  we know that a minimal set of generators of $\REQ$ is given by 
	$$
	\big\{ y^2T_2-x^2T_3,\quad y^3T_1-x^3T_2,\quad xT_2^2-yT_1T_3,\quad
	yT_2^3-xT_1T_3^2,\quad
	T_2^5-T_1^2T_3^3
	\big\},
	$$
	so a minimal set of generators for $\EEQ$ is given by
	$$
	\big\{xT_2^2-yT_1T_3,\quad
	yT_2^3-xT_1T_3^2,\quad
	T_2^5-T_1^2T_3^3
	\big\},
	$$
	We make the following session in Macaulay2:
	\vspace*{.4cm}
	\begin{spacing}{0.75}
		{   \normalfont
			\fontfamily{cmtt}
			\selectfont
			\hspace*{-.4cm}i1 : \verb|R = QQ[x,y]|\\
			o1 = \hspace{.5pt} \verb|R|\\
			o1 : \verb|PolynomialRing|\\
			i2 : \verb|load "bFunctionRees.m2"|\\
			i3 : \verb|I = ideal(x^5, x^2*y^3, y^5)|\\
			\hspace*{2.35cm}5 \hspace*{.14cm} 2 3 \hspace*{1pt} 5\\
			o3 = \hspace{.5pt} ideal(x\;, x y, y )\\
			o3 : \verb|Ideal of R|\\
			i4 : \verb|for p from 2 to 5 do << factorBFunction bFunctionRees(I, p) << endl;|\\
			(s)(s + 1)(s + 2)\\
			(s)(s + 1)(s + 2)\\
			(s)(s + 1)(s + 2)\\	
			(s)(s + 1)(s + 2)(s + 3)\\
		}	
	\end{spacing}	

	From \autoref{expression_bFunction} we see that for $p=2,\ldots,4$ we have $\EEQ_{p,q}\neq 0$ if and only if $1 \le q \le 3$, and that $\EEQ_{5,q} \neq 0$ if and only if $0 \le q \le 3$.
\end{example}

\begin{example}
	We assume that in \autoref{data} we have $\mu=1$. 
	In this case it is known (see e.g. \cite[Theorem 2.3]{COX_HOFFMAN_WANG} or \cite[Proposition 3.1]{LAURENT}) that the minimal generators of $\REQ$ have bidegrees
	$$
	(1,1),\;(1,d-1),\;(2,d-2),\;(3,d-3),\;\ldots,\;(d,0).
	$$
   We can make an interesting session with ideals of this form created randomly, we take the particular case $\mu=1$ and $d=7$:
   \vspace*{.4cm}
	\begin{spacing}{0.75}
		{   \normalfont
			\fontfamily{cmtt}
			\selectfont
	\hspace*{-.6cm} i1 : \verb|R = QQ[x,y]|\\
	o1 = \hspace{.5pt} \verb|R|\\
	o1 : \verb|PolynomialRing|\\
	i2 : \verb|load "bFunctionRees.m2"|\\
	i3 : \verb|A = matrix{{random(1,R),random(6,R)},{random(1,R),random(6,R)},|\\
	\hspace*{1.5cm}\verb|{random(1,R),random(6,R)}};|\\
	\hspace*{2.6cm}3 \hspace*{.7cm} 2\\
	o3 : Matrix R \verb|<---| R\\
	i4 : \verb|I =  minors(2, A);|\\
	o4 : Ideal of R\\
	i5 : assert(codim I == 2);\\
	i6 : \verb|for p from 2 to 7 do << factorBFunction bFunctionRees(I, p) << endl;|\\
	(s)\\
	(s)(s + 1)\\
	(s)(s + 1)(s + 2)\\
	(s)(s + 1)(s + 2)(s + 3)\\
	(s)(s + 1)(s + 2)(s + 3)(s + 4)\\
	(s)(s + 1)(s + 2)(s + 3)(s + 4)(s + 5)\\			
		}	
	\end{spacing}
Here we need to check ({\normalfont\fontfamily{cmtt}
	\selectfont assert(codim I == 2);}) that the created ideal $I$ has height $2$, although it is extremely improbable that this is not the case.
\end{example}
	
\section*{Acknowledgments}
I am grateful to my PhD advisor, Carlos D'Andrea, for his support and guidance, for suggesting the problem, and for the many helpful discussions.
I am grateful to Laurent Bus\'e and Santiago Zarzuela, for their patience,  and for the helpful discussions and suggestions.
I am  thankful to Francisco Jes\'us Castro Jim\'enez for several suggestions and improvements on an early draft of this paper.

I thank Josep \`Alvarez Montaner, Alessio Caminata and Ricardo Garc\'ia for useful discussions. 
The use of \textit{Macaulay2} \cite{M2} was a driving force in the preparation of this paper.

\bibliographystyle{elsarticle-num} 

%\bibliography{references}

\end{document}